\documentclass{amsart}
\usepackage{amssymb}
\usepackage{amscd}
\usepackage{amsmath}
\usepackage{enumerate}
\usepackage[dvips]{graphicx}
\newtheorem{theo}{Theorem}
\newtheorem{lema}[theo]{Lemma}
\newtheorem{cor}[theo]{Corollary}
\newtheorem{prop}[theo]{Proposition}

\newtheorem{definition}[theo]{Definition}

\newcommand{\CC}{{\mathbb{C}}}

\newcommand{\QQ}{{\mathbb{Q}}}
\newcommand{\RR}{{\mathbb{R}}}
\newcommand{\KK}{{\mathbb{K}}}
\newcommand{\SSS}{{\mathbb{S}}}
\newcommand{\ZZ}{{\mathbb{Z}}}

\newcommand{\calU}{{\mathcal{U}}}

\newcommand{\comp}{{\circ}}

\begin{document}
\title[Nash problem for surfaces]{Nash problem for surfaces}
\author{Javier Fern\'andez de Bobadilla}
\address{ICMAT. CSIC-UAM-UCM-UC3M. Departamento de Algebra. Facultad de Ciencias Matem\'aticas. Universidad Complutense de Madrid. Plaza de Ciencias 3. 28040 Madrid}
\email{javier@mat.csic.es}
\author{Mar\'ia Pe Pereira}
\address{ICMAT. CSIC-UAM-UCM-UC3M}
\curraddr{Institut de Math\'ematiques de Jussieu\\
\'Equipe G\'eom\'etrie et Dynamique} 
\email{maria.pe@mat.ucm.es}
\thanks{Second author is supported by Caja Madrid. Research partially supported by the ERC Starting Grant project TGASS and by Spanish Contracts MTM2007-67908-C02-02 and MICINN2010-2170-C02-01. The authors are grateful to the Faculty of Ciencias Matem\'aticas of Universidad Complutense de Madrid for excellent working conditions.}
\date{14-2-2011}
\begin{abstract}
We prove that Nash mapping is bijective for any surface defined over an algebraically closed field of
characteristic $0$.
\end{abstract}


\maketitle

\section{Introduction}

Nash problem~\cite{Na} was formulated in the sixties (but published later) in the attempt to understand the relation between the 
structure of resolution of singularities of an algebraic variety $X$ over a field of characteristic $0$ and the space of arcs (germs of parametrized curves) in the variety.
He proved that the space of arcs centred at the singular locus (endowed with an infinite-dimensional algebraic variety structure) 
has finitely many irreducible components and proposed to study the relation of these components with the 
essential irreducible components of the exceptional set a resolution of singularities. 

An irreducible component $E_i$ of the exceptional divisor of a
resolution of singularities is called essential, if given any other resolution 
the birational transform of $E_i$ to the second resolution is an irreducible component of the exceptional divisor. Nash defined a mapping from the set of irreducible components of the
space of arcs centred at the singular locus to the set of essential components of a resolution as follows: he assigns to each component $W$ 
of the space of arcs centred at the singular locus the unique component of the exceptional divisor
which meets the lifting of a generic arc of $W$ to the resolution. Nash established the injectivity of this mapping. For the case of surfaces it seemed possible for him that the mapping is also surjective, and posed
the problem as an open question. He also proposed to study the mapping in the higher dimensional case.
Nash resolved the question positively for the $A_k$ singularities. As a general reference for Nash problem the reader may look at~\cite{Na} and~\cite{IK}.

Besides Nash problem, the study of arc spaces is interesting because it lays the
foundations for motivic integration and because the study of its geometric properties reveals properties of the underlying varieties (see papers of Denef, Loeser, de Fernex, Ein, Ishii, Lazarsfeld, Mustata, Yasuda and others).

It is well known that birational geometry of surfaces is much simpler than in higher dimension. This fact reflects on Nash problem: Ishii and Kollar showed in~\cite{IK} a 4-dimensional example with non-bijective Nash mapping. In the same paper they showed the bijectivity of the Nash mapping for toric singularities of 
arbitrary dimension. Other advances in the higher dimensional case include~\cite{PP2},~\cite{Go},~\cite{LR}.

On the other hand bijectivity of the Nash mapping has been shown for many classes of surfaces (see \cite{Go},\cite{IK},\cite{I1},\cite{I2},\cite{Le},\cite{LR1},\cite{Mo},\cite{Pe},\cite{Pl},\cite{Pl2},\cite{PlSp}, \cite{PP1},\cite{Re1},\cite{Re2}). The techniques leading to the proof of each of these cases are 
different in nature, and the proofs are often complicated. It is worth to notice that even for the case of
the rational double points not solved by Nash a complete proof has to be awaited
until last year: see~\cite{Pl2},
~\cite{PlSp} and~\cite{Pe}; in the last paper the result is proved for all quotient surface singularities.

In this paper we resolve Nash question for surfaces:\\

\noindent\textbf{Main Theorem.}
{\em Nash mapping is bijective for any surface defined over an algebraically closed field of characteristic $0$.}\\

The core of the result is the case of normal surface singularities. After settling this case we deduce from it the general surface
case following a suggestion by C. Pl\'enat and M. Spivakovsky.

The proof is based on the use of convergent wedges and topological methods. A wedge is a uniparametric family of arcs. The use of wedges in connection 
to Nash problem was proposed by M. Lejeune-Jalabert~\cite{Le}. Later A. Reguera~\cite{Re}, 
building onto the fundamental Lemma of motivic integration by J. Denef and F. Loeser~\cite{DL}, proved a 
characterization of components which are at the image of the Nash map in terms of formal wedges defined over
fields which are of infinite transcendence degree over the base field.
In~\cite{Bo} it is proved a characterization of the image of the Nash mapping for surfaces
in terms of convergent (or even algebraic) wedges defined over the base field, which is the starting point of this article. 
In the same paper it is shown 
that Nash problem is of topological nature. Independently, in~\cite{LR}, it is given a different sufficient condition for being at the image of the Nash map in terms of formal wedges defined over the base field; this
condition holds in arbitrary dimension, but it is weaker than the one given in~\cite{Bo} for the surface 
case. In~\cite{Pe} the second author settles Nash question  for quotient surface singularities. The present paper is inspired in the ideas of~\cite{Pe}, more concretely in the use of representatives of wedges and in the use of topological methods.

The proof in~\cite{Pe} uncovers a nice phenomenon of the deformation theory of curves
in surfaces, which is worth to be studied on its own, and which perhaps could give
a different proof for the bijectivity of Nash mapping for surfaces.

The idea of our proof is as follows: let $(X,O)$ be a normal surface singularity and
$$\pi:\tilde{X}\to (X,O)$$
be the minimal resolution of singularities.
By a Theorem of~\cite{Bo} if Nash mapping of $(X,O)$ is not bijective there exists a convergent wedge
$$\alpha:(\CC^2,O)\to (X,O)$$
with certain precise properties. As in~\cite{Pe}, taking a suitable representative we may view $\alpha$ as a
uniparametric family of mappings
$$\alpha_s:\calU_s\to (X,O)$$
 from a family of domains $\calU_s$ to $X$ with the property that each $\calU_s$ is diffeomorphic to a disc.
For any $s$ we consider the lifting
$$\tilde{\alpha}_s:\calU_s\to\tilde{X}$$
to the resolution. Notice that $\tilde{\alpha}_s$ is the normalization mapping of the image curve. 

On the other hand, if we denote by $Y_s$ the image of $\tilde{\alpha}_s$ for $s\neq 0$, then we may consider the
limit divisor $Y_0$ in $\tilde{X}$ when $s$ approaches $0$. This limit divisor consists of the union of the image of 
$\tilde{\alpha}_0$ and certain components of the exceptional divisor of the resolution whose multiplicities are easy to be computed. We prove an upper
bound for the Euler characteristic of the normalization of any reduced deformation of $Y_0$ in terms of the following data:
the topology of $Y_0$, the multiplicities of its components and
the set of intersection points of $Y_0$ with the generic member $Y_s$ of the deformation. Using this bound we show that the Euler characteristic of the normalization of $Y_s$
is strictly smaller than one. This contradicts the fact that the normalization is a disc.

In the last Section we deduce the general case from the normal case following a remark by C. Pl\'enat and M. Spivakovsky.
 
\section{Preliminaries}
\subsection{}\label{res} 
Let $(X,O)$ be a complex analytic normal surface singularity. Let 
$$\pi:(\tilde{X},E)\to (X,O)$$
be the minimal resolution of singularities, which is an isomorphism outside the exceptional divisor $E:=\pi^{-1}(O)$. Consider the decomposition 
$E=\cup_{i=0}^r E_i$ of $E$ into irreducible components. These irreducible components are the essential components of $(X,O)$.

The germ $(X,O)$ is embedded in an ambient space $\CC^N$. Denote by $B_{\epsilon}$ the closed ball of radius $\epsilon$ centred at the 
origin and by $\SSS_{\epsilon}$ its boundary sphere. Take a Milnor radius $\epsilon_0$ for $(X,O)$ in $\CC^N$ (in particular $X\cap B_{\epsilon_0}$ has conical differentiable structure). From now on we will denote by $X$ a representative $X\cap B_{\epsilon_0}$ of $X$ and by $\tilde{X}$ the resolution of singularities $\pi^{-1}(X)$.  
In these conditions 
the space $\tilde{X}$ admits the exceptional divisor $E$ as a deformation retract. Hence the 
homology group $H_2(\tilde{X},\ZZ)$ is free and generated by the classes of the irreducible components $E_i$. Since $\tilde{X}$ is a 
smooth $4$-manifold there is a symmetric intersection product
$$\centerdot:H_2(\tilde{X},\ZZ)\times H_2(\tilde{X},\ZZ)\to\ZZ.$$
The intersection product is negative definite since it is the intersection product of a resolution of a surface singularity. 

\subsection{}\label{wedgeformulation} We recall some terminology and results from~\cite{Bo}. Consider coordinates $(t,s)$ in the germ $(\CC^2,O)$. A \emph{convergent wedge} is a complex analytic germ
$$\alpha:(\CC^2,O)\to (X,O)$$
which sends the line $V(t)$ to the origin $O$. Given a wedge $\alpha$ and a parameter value $s$, the arc
$$\alpha_s:(\CC,0)\to (X,O)$$
is defined by $\alpha_s(t)=\alpha(t,s)$. The arc $\alpha_0$ is called {\em the special arc} of the wedge. For small enough $s\neq 0$
the arcs $\alpha_s$ are called \emph{generic arcs} of $\alpha$.

Any arc 
$$\gamma:(\CC,0)\to (X,O)$$
admits a unique lifting $\tilde{\gamma}$ to $(\tilde{X},O)$.

\begin{definition}[\cite{Bo}]
\label{wedadj}
A convergent wedge $\alpha$ {\em realizes an adjacency} from $E_j$ to $E_i$ (with $j\neq i$) if and only if the lifting $\tilde{\alpha}_0$ of the special arc meets $E_i$ transversely at a non-singular point of $E$ and the lifting $\tilde{\alpha}_s$ of a generic arc satisfies $\tilde{\alpha}_s(0)\in E_j$.
\end{definition}

Our proof is based in the following Theorem, which is the implication ``~$(1)~\Rightarrow~(a)$~'' of Corollary~B~of~\cite{Bo}:
\begin{theo}[\cite{Bo}]
\label{nashwedges}
An essential divisor $E_i$ is in the image of the Nash mapping if there is no other essential divisor 
$E_j\neq E_i$ such that there exists a convergent wedge realizing an adjacency from $E_j$ to $E_i$.
\end{theo}

\subsection{}\label{lefschetz} The previous theorem allows to address Nash question in the complex analytic case. 
Supose that $(X,O)$ is a singularity of a normal algebraic surface defined over an algebraically closed field $\KK$ 
of characteristic $0$. It is well known that $(X,O)$ may be defined over a field $\KK_1$ which is a finite 
extension of $\QQ$, and hence admits an embedding into $\CC$. Let $\bar{\KK}_1$ the algebraic closure of 
$\KK$. We have then two field embeddings $\overline{\KK}_1\subset\KK$ and $\overline{\KK}_1\subset\CC$.
Since the bijectivity of the Nash mapping does not change by extension of algebraically closed fields we 
deduce that if we prove the bijectivity of the Nash mapping for any complex analytic normal surface 
singularity, then it holds for any normal surface singularity defined over a field of characteristic equal to $0$.

\subsection{}\label{representatives}
Following~\cite{Pe} we shall work with representatives rather than germs in order to get richer information about the geometry of the possible wedges. 
Remember that $X$ stands for a Milnor representative $X\cap B_{{\epsilon}_0}$ for a Milnor radius $\epsilon_0$ for $(X,O)$. 

Given any non-constant arc germ $\gamma:(\CC,O)\to (X,O)$ there exist a 
suitable representative 
$$\gamma:\Omega\to X,$$
where $\Omega$ is an open neighbourhood of the origin of $\CC$,
and a positive radius $\epsilon$ such that the mapping 
$$\gamma:\gamma^{-1}(B_{\epsilon})\to B_{\epsilon}$$ 
is proper and transversal to the sphere $\SSS_\rho$ for any $0<\rho\leq\epsilon$ and, moreover, the preimage $\gamma^{-1}(B_\epsilon)$ is connected. 
In this case $\gamma^{-1}(B_\epsilon)$ is diffeomorphic to a closed disc.

\begin{definition}[\cite{Pe}]
A {\em Milnor representative} of $\gamma$ is a representative of the form 
$$\gamma|_{D}:D\to X$$
where $D$ is the disc $\gamma^{-1}(B_\epsilon)$ for $\epsilon$ as above.
\end{definition}

Denote by $D_\delta$ the disc of radius $\delta$ centred at the origin of $\CC$.

Given a wedge $\alpha$ with non-constant special arc $\alpha_0$, we consider the radius $\epsilon$ associated to a Milnor representative of $\alpha_0$
such that $\alpha$ is defined in a neighbourhood of $D\times D_\delta$ for a positive $\delta$.
The radius $\delta$ can be chosen small enough
so that the mapping $\alpha_s$ is transversal to $\SSS_\epsilon$ for
any $s\in D_\delta$. We consider the mapping
$$\beta:(\CC^2,O)\to(\CC^N\times D_\delta,(O,0))$$
given by $\beta(t,s):=(\alpha(t,s),s)$ and define $\calU_{\epsilon,\delta}:=\beta^{-1}(B_\epsilon\times D_\delta)$. We have that the 
restriction to $\calU_{\epsilon,\delta}$ of the natural projection to $D_\delta$ is a smooth trivial fibration with fibre diffeomorphic
to a closed disc.

\begin{definition}[\cite{Pe}]
\label{goodrepresentative}
A {\em good representative} of a wedge $\alpha$ is a representative of the form 
$$\alpha|_{\calU_{\epsilon,\delta}}:\calU_{\epsilon,\delta}\to X$$ 
with 
$\calU_{\epsilon,\delta}$ and $\delta$ chosen as in the previous discussion.
\end{definition}

From now on we only deal with wedges with non-constant special arc, which, therefore, have good representatives. 

Given a good representative $\alpha|_{\calU_{\epsilon,\delta}}$, we denote the corresponding $\calU_{\epsilon,\delta}$ simply by $\calU$ and consider the representative 
$$\beta|_{\calU}:\calU\to X\times D_\delta.$$ 
  Given any $s\in D_\delta$ we denote by
$\calU_s$ the fibre by the natural projection of $\calU$ onto $D_\delta$; it is a region in $\CC$ diffeomorphic to a disc. 
The fact that $\calU_s$ is a disc is a key in the proof as it was in the final step of the proof of the main result of~\cite{Pe}. 
\section{Wedges and divisors}

%
\subsection{}\label{familiadivisores}
Given a wedge $\alpha$ realizing an adjacency, consider a good representative $\alpha|_{\calU}$ as in Definition~\ref{goodrepresentative}. 

Since $\beta|_{\calU}$ is proper, its image $H:=\beta(\calU)$ is a $2$-dimensional closed 
analytic subset of $X\times D_\delta$. For any $s\in D_\delta$ the fibre $H_s$ by the natural projection onto $D_\delta$ is the image of the representative
$$\alpha_s|_{\calU_s}:\calU_s\to X.$$

Given the minimal resolution of singularities 
$$\pi:\tilde{X}\to X$$ we consider the mapping
$$\sigma:\tilde{X}\times D_\delta\to X\times D_\delta.$$
The inverse image $\sigma^{-1}(H)$ is a divisor in the smooth $3$-fold $\tilde{X}\times D_\delta$, which can be decomposed as
$$\sigma^{-1}(H)=Y+\sum_{i=0}^r n_i(E_i\times D_\delta),$$
where $Y$ denotes the strict transform of $H$.

For any $s\in D_{\delta}$ there is a unique lifting 
\begin{equation}
\tilde{\alpha}_s:\calU_s\to Y_s\subset\tilde{X}
\end{equation}
of $\alpha_s|_{\calU_s}$ to $\tilde{X}$. 

For $s\neq 0$ the fibre $Y_s$ is reduced and coincides with the image of the mapping $\tilde{\alpha}_s$. Therefore 
$\tilde{\alpha}_s$ is the normalization mapping of the curve $Y_s$. For $s=0$ the divisor $Y_0\subset X$ decomposes as
\begin{equation}
\label{decomposition}
Y_0=Z_0+\sum_{i=0}^ra_iE_i,
\end{equation}
where $Z_0$ is the image of the lifting $\tilde{\alpha_0}$.

\subsection{}\label{sec:sistemalineal} If $\alpha$ is a wedge realizing the adjacency from $E_j$ to $E_0$ with $j\neq 0$, then by definition, the lifting $\tilde{\alpha}_0$ meets $E_0$ transversely. In particular, $Z_0\centerdot E_0=1$ and 
$Z_0\centerdot E_i=0$ for $i>0$, where $Z_0$ is as in Formula~(\ref{decomposition}).

Since the divisor $Y_s$ is a deformation of the divisor $Y_0$ we have the equality 
\begin{equation}
\label{eqns}Y_0\centerdot E_i=Y_s\centerdot E_i\end{equation}
for any $i$. Denote by $b_i$ the intersection product of $Y_s\centerdot E_i$ and by $M$ the intersection matrix of the intersection product in $\tilde{X}$. Then, (\ref{eqns}) can be expressed as follows:
\begin{equation}
\label{eq:sistemalineal}
M(a_0,...,a_r)^t=(-1+b_0,b_1,...,b_r)^t.
\end{equation}

In the terminology of~\cite{Pe} the number $b_i$ is the number
of \emph{returns} of the wedge \emph{through the divisor} $E_i$: it is the number of points $p\in\alpha_s|_{\calU_s}^{-1}(O)$ for which the lifting to $\tilde{X}$ of the germ at $p$ of $\alpha_s|_{\calU_s}$ meets $E_i$. 
The use of returns was one of key ideas introduced by the second author in \cite{Pe}. 

Since $\alpha$ realizes an adjacency from $E_j$ to $E_0$ we have more restrictions about $b_i$'s and $a_i$'s. They can be seen as consequences of the following lemma: 
\begin{lema}
\label{entradasnegativas}
All the entries of the inverse matrix $M^{-1}$ are non-positive.
\end{lema}
\begin{proof}
The matrix $-M$ is symmetric, positive definite, and such that any non-diagonal entry is non-positive. Hence, if endow $\RR^r$ with the standard euclidean product then there is a basis ${v_1,...,v_r}$ such that the angle formed by any two
different vectors of the base is at least $\pi/2$, and the matrix $-M$ is the matrix of scalar products of pairs of vectors of the
basis. Therefore the inverse matrix $-M^{-1}$ is the matrix of scalar products of pairs of vectors of a basis of vectors such that
the angle formed by any two of the vectors is at most $\pi/2$. This implies that all the entries of $-M^{-1}$ are non-negative.
\end{proof}
Hence, if we require in (\ref{eq:sistemalineal}) that each $b_i$ and each $a_i$ are non-negative integers, then we get that $b_0$ has to be equal to $0$ or to $1$, and in this last case we get that $b_1=...=b_r=0$. 
Hence, we have the following immediate consequence:
\begin{cor}
\label{a_0}
If $\alpha$ is a wedge realizing an adjacency from $E_j$ to $E_0$ (with $j\neq 0$), and $(b_0,...,b_r)$ are the intersection numbers $Y_s\cdot E_i$ associated with the generic member of a good wedge representative as in~(\ref{eq:sistemalineal}), then $b_0$ is equal to $0$. Moreover $a_0$ is positive, that is the divisor $E_0$ appears in the support of $Y_0$.
\end{cor}
\begin{proof}
Since $\alpha$ realizes an adjacency from $E_j$ to $E_0$ we have $b_j\neq 0$. Then $b_0=0$. Now in the first row of system~(\ref{eq:sistemalineal}) in order to have the equality $b_0=0$ we need that $\sum_{j=0}^r a_{j}k_{0,j}=-1$. By definition all $a_j$ and all $k_{0,j}$ except $k_{0,0}$ are non-negative. This implies that $a_0$ is different from $0$. 
\end{proof}

\subsection{}\label{firstrestrictions}
The equality~(\ref{eq:sistemalineal}) can be viewed as a linear system whose indeterminates are $a_0,...,a_r$. It can be used to prove that wedges realizing certain adjacencies with certain prescribed returns do not exist (we are using the terminology of~\cite{Pe}).
The method is as follows: the adjacencies and the prescribed returns determine $b_0,b_1,...,b_r$. The existence of the wedge is 
impossible if the solution of the linear system has either a negative or a non-integral entry.

Using this method it is possible to prove the bijectivity of Nash mapping for many singularities (toric, dihedral...), but it does not
suffice for all of them. It is interesting to compare this method with the methods of~\cite{Pe} for the $E_8$ singularity. The set of 
adjacencies with prescribed returns which this method is not able to rule coincide precisely with the list of 25 adjacencies with 
with prescribed returns that the second author is not able to rule out only with intersection multiplicity methods. 

\section{Euler characteristic estimates}\label{sec:eulerchar}

Let $\tilde{X}$ be a compact domain with smooth boundary in a compact complex surface. Let 
$$Y_0=\sum_{i=0}^mc_iZ_i+\sum_{i=0}^ra_iE_i$$
be a divisor in $\tilde{X}$, where the $E_i$'s are compact prime divisors contained in the interior of $\tilde{X}$, and the $Z_i$'s
are prime divisors meeting transversely the boundary of $\tilde{X}$. We consider a deformation $Y_s$ of the divisor $Y_0$ such that 
$Y_s$ is reduced and transversal to the boundary $\partial\tilde{X}$ for $s\neq 0$. Let $$n:\calU_s\to Y_s$$ be the normalization of $Y_s$. 
In this section we bound the Euler characteristic of the normalization $\calU_s$ in terms 
of the topology of the reduced divisor associated with $Y_0$, the multiplicities $c_i$ and $a_i$ and the number of intersection 
points of $Y_0$ with $Y_s$, for $s\neq 0$.

We do first the case when $Y_0$ is a normal crossing divisor. We denote by $(Y_0)^{red}$ the reduced divisor associated with $Y_0$.

\subsection{Local normal crossings case}\label{normalcrossingslocal}
In this case $\tilde{X}$ is a ball $B_\epsilon$ 
centred at the origin of $\CC^2$, and $Y_0$ is defined by $f_0=x^ay^b=0$, where $x$ and $y$ are the coordinates of $\CC^2$. The divisor $Y_s$ is defined by $f_s=0$, where $f_s$ is a $1$-parameter holomorphic deformation of $f_0$
such that $f_s$ is reduced for $s\neq 0$.
We have the following bound:   
\begin{lema}
If $s$ is small enough then the Euler characteristic of the normalization $\calU_s$ of $Y_s$ satisfies:
\begin{equation}
\chi(\calU_s)\leq\sum_{p\in Y_s\cap Y_0}I_p((Y_0)^{red},Y_s).
\end{equation}
\end{lema}
\begin{proof}
The only connected orientable surface with boundary which has positive Euler characteristic is the disc. Hence $\chi(\calU_s)$ is bounded
above by the number of connected components of $\calU_s$ which are discs. 

Let $W_s$ be an irreducible component of $Y_s$ whose normalization is a disc. Its boundary $W_s\cap\SSS_\epsilon$
is a circle which deforms to one of the components of $Y_0\cap \SSS_\epsilon$, that is either to $V(x)\cap\SSS_\epsilon$ or
to $V(y)\cap\SSS_\epsilon$. Both cases are symmetric. In the first case the equation $g_s$ of $W_s$ degenerates to $x^c$ for a certain
$c\leq a$, that is $g_0=x^c$. Thus the circle $W_s\cap\SSS_\epsilon$ loops $c$ times around the $V(y)$, and hence represents
a non-trivial element in $\pi_1(B_\epsilon\setminus V(y))$. The normalization of the component $W_s$ is a mapping from a disc to 
$W_s$. If $W_s$ does not meet $V(y)$ the circle $W_s\cap\SSS_\epsilon$ would be a trivial element in 
$\pi_1(B_\epsilon\setminus V(y))$, and this is not the case.

We conclude that each component of $Y_s$ whose normalization is a disc has at least one intersection point with the union of the axis. 
This proves the lemma.
\end{proof}

\subsection{Global normal crossings case}\label{normalcrossingsglobal}
We assume $Y_0$ to be a normal crossings divisor. Define 
$$\dot{E}_i=E_i\setminus {\rm Sing}((Y_0)^{red}),$$
$$\dot{Z}_i=Z_i\setminus {\rm Sing}((Y_0)^{red}),$$
for any $i$.
\begin{lema}
\label{cotanormalcrossingsglobal}
If $s$ is small enough then the Euler characteristic of the normalization $\calU_s$ of $Y_s$ satisfies
\begin{equation}\label{eq:normalcrossingsglobal}
\chi(\calU_s)\leq\sum_{i=0}^mc_i\chi(\dot{Z}_i)+\sum_{i=0}^ra_i\chi(\dot{E}_i)+\sum_{p\in Y_s\cap Y_0}I_p((Y_0)^{red},Y_s).
\end{equation}
\end{lema}
\begin{proof}
Choose small balls $B_1,...,B_k$ inside $\tilde{X}$ centred in each of the singular points of $(Y_0)^{red}$. Choose tubular 
neighbourhoods $T_i$ (respectively $T'_i$) around each component $Z_i$ (respectively $E_i$), which are so small that their boundaries
meet the boundary of each of the balls $B_j$ transversely. Define
$$U_i:=T_i\setminus\cup_{j=1}^k B_j,$$
$$W_i:=T'_i\setminus\cup_{j=1}^k B_j.$$
We have product structures 
$$U_i\cong (Z_i\setminus\cup_{j=1}^k B_j)\times D$$
and 
$$W_i\cong (E_i\setminus\cup_{j=1}^k B_j)\times D.$$
Let 
$$\rho_i:U_i\to Z_i\setminus\cup_{j=1}^k B_j,$$ 
$$\kappa_i:W_i\to E_i\setminus\cup_{j=1}^k B_j$$
be the projections onto the first factors. 

If $s$ is small enough the compositions
$$\rho_i\comp n|_{n^{-1}(Y_t\cap U_i)}n^{-1}(Y_t\cap U_i)\to Z_i\setminus\cup_{j=1}^k B_j$$
$$\kappa_i\comp n|_{n^{-1}(Y_t\cap W_i)}n^{-1}(Y_t\cap U_i)\to E_i\setminus\cup_{j=1}^k B_j$$
are branched covers of degree $c_i$ and $a_i$ for any $i$. By Hurwitz formula we find
$$\chi(n^{-1}(Y_t\cap U_i))\leq c_i\chi(\dot{Z}_i)$$
and
$$\chi(n^{-1}(Y_t\cap W_i))\leq a_i\chi(\dot{E}_i).$$

Since  $n^{-1}(Y_s\cap\partial B_i)$ is a union of circles for any $i$ and the Euler characteristic of a circle is equal to $0$, 
we have 
$$\chi(\calU_s)=\sum_{i=1}^k\chi(n^{-1}(Y_t\cap B_i))+\sum_{i=0}^m\chi(n^{-1}(Y_t\cap U_i))+\sum_{i=0}^r\chi(n^{-1}(Y_t\cap W_i)).$$
Using the local bound obtained in paragraph ~\ref{normalcrossingslocal} and the bounds above we get the required bound. 
\end{proof}

\subsection{General local case}\label{generallocal}
In this case $Y_0$ is defined by $f_0=\prod_{i=0}^mg_i^{c_i}=0$, where the $g_i$ are irreducible and reduced analytic function germs. 
We denote by $\mu_i$ the Milnor number of $g_i$ at the origin. 
We take a Milnor ball $B_\epsilon$ for $f_0$ as the space $\tilde{X}$. The divisor $Y_s$ is defined by $f_s=0$, where $f_s$ is a $1$-parameter holomorphic deformation of $f_s$ such that $f_s$ is reduced for $s\neq 0$. We consider a sufficiently small $\delta$ so that 
$f_0^{-1}(\delta)\cap B_\epsilon$ is the Milnor fibre of $f_0$. We will use and generalize in certain sense the following equality that was proved in~\cite{Me}. We start by giving an alternative proof of the equality.
\begin{lema}[\cite{Me}]
\label{cotafibramilnor}
The Euler characteristic of the Milnor fibre of $f_0$ is equal to:
\begin{equation}\label{eq:milnor_char1}\chi(f_0^{-1}(\delta)\cap B_\epsilon)=\sum_{i=0}^mc_i(1-\mu_i-I_O(g_i,\prod_{j\neq i}g_j)).\end{equation}
\end{lema}
\begin{proof}
Given a vector $v$ of $\CC^2$ we denote by $\tau_v$ the translation of $\CC^2$ associated with $v$. 
We choose $m$ vectors $v_1,...,v_m$ in $\CC^2$ such that for any $t$ small enough and $i\neq j$ the curves $V(g_i\comp\tau_{tv_i}-t)$
and $V(\prod_{j\neq i}g_j\comp\tau_{tv_j}-t)$ meet transversely in $B_\epsilon$. 

Consider the deformation $F_t:=\prod_{i=0}^m(g_i\comp\tau_{tv_i}-t)^{c_i}$. An easy local argument shows that
for small enough $t$ and any $s\in D_\delta\setminus\{0\}$ the set $F_t^{-1}(s)$ is smooth at the meeting points 
with $\partial B_{\epsilon}$ and transverse to it. This implies the existence of a finite subset of critical values 
$\Delta_t$ of $D_\delta$ such that the restriction
$$F_t:B_\epsilon\cap F_t^{-1}(D_\delta\setminus\Delta_t)\to D_\delta\setminus\Delta_t$$
is a locally trivial fibration with fibre diffeomorphic to the Milnor fibre of $f_0$.
See Theorem~2.2~of~\cite{Bo0} for a proof of these facts in a much more general context.

Fix a small enough $t$ different from $0$. We view $F_t^{-1}(s)$ as a deformation of the normal crossings divisor $F_t^{-1}(0)$ inside 
$B_\epsilon$ and study it like in the global normal crossings case.
The irreducible components of this divisor are $Z_i=V(g_i\comp\tau_{tv_i}-t)$, for $i=0,...,m$. The component $Z_i$ is a translation of the the Milnor fibre of $g_i$, and, hence, its Euler characteristic is equal to $1-\mu_i$. Consequently, using that 
the curve $Z_i$ meets transversely the union $\cup_{j\neq i} Z_i$ and the conservativity of intersection multiplicity, we obtain 
$$\chi(\dot{Z}_i)=1-\mu_i-I_O(g_i,\prod_{j\neq i}g_j).$$

Observe that the piece of the Milnor fibre contained in a neighbourhood of a singularity of $F_t^{-1}(0)$ is a union of
cylinders because locally $F_t^{-1}(0)$ is normal crossings. Decomposing the Milnor fibre as in
Lemma~\ref{cotanormalcrossingsglobal} and adding the Euler characteristics of the corresponding pieces 
we obtain the equality.
\end{proof}

After this Lemma we can prove the Euler characteristic bound that we want:
\begin{prop}
\label{cotagenerallocal}
If $s$ is small enough we have
\begin{equation}\label{eq:cotagenerallocal}
\chi(\calU_s)\leq\sum_{i=0}^mc_i(1-\mu_i-I_O(g_i,\prod_{j\neq i}g_j))+\sum_{p\in Y_s\cap Y_0}I_p((Y_0)^{red},Y_s).
\end{equation}
\end{prop}
\begin{proof}$\quad$

\textsc{A particular case:} the divisor $Y_s$ does not meet the origin for $s\neq 0$. 

In order to reduce the problem to the global normal crossings case we consider the minimal embedded resolution 
$$\pi:\tilde{X}\to B_\epsilon$$ 
of $V(f_0)$. Let $\{E_i\}_{i=1}^{r}$ be the irreducible components of the exceptional divisor. For any $s\in D_\delta$ we denote by $V_s$ the pullback of $Y_s$ by $\pi$. 
Since the divisor $Y_s$ does not meet the origin when $s\neq 0$ we have that it is isomorphic to $V_s$ and that $V_s$ does not meet the exceptional divisor of $\pi$. Then it is enough to prove the bound for the Euler characteristic of the divisor 
$V_s$ for $s\neq 0$.

The divisor $V_0$ decomposes as $V_0=\sum_{i=0}^m c_iZ_i+\sum_{i=0}^ra_iE_i$, where the $c_i$'s depend only on  the resolution of singularities of $f_0^{-1}(0)$ and the $a_i$'s are deduced from the $c_i$'s solving the linear system derived from the identities $V_0\cdot E_i=0$ (notice that $V_s$ does not meet any $E_i$ because $Y_s$ does not meet 
the origin and then $V_s\cdot E_i=0$ for all $i$). 

Using the bound obtained in paragraph~\ref{normalcrossingslocal} and the fact that $\dot{Z}_i$ is a punctured disc for any $i$ we obtain
\begin{equation}\label{aux}\chi(\calU_s)\leq\sum_{i=0}^ra_i\chi(\dot{E}_i)+\sum_{p\in V_0\cap V_s}I_p(V_s,(V_0)^{red}).\end{equation}
Using the fact that the number of intersection points of $Y_s$ and $(Y_0)^{red}$
counted with multiplicity coincides with the number of intersection
points of $V_s$ and $(V_0)^{red}$ counted with multiplicity, after Lemma~\ref{cotafibramilnor}, in order to prove the proposition it only rest to check that the first sum of the right side of~(\ref{aux}) coincides with the Euler characteristic of the Milnor fibre of $f_0$. 

For this we observe that the divisor $V_0=\sum_{i=0}^m c_iZ_i+\sum_{i=0}^ra_iE_i$ is equal to the total transform of $V(f_0)$ by
the modification $\pi$. This is because the coefficients $a_i$ are also characterized by the 
equalities $V_0\cdot E_i=0$ for any $i$. The Euler characteristic of the Milnor fibre is given then by 
\begin{equation}\label{eq:milnor_char2}\chi(f_0^{-1}(s))=\sum_{i=0}^ra_i\chi(\dot{E}_i).\end{equation}
Indeed, if $W_s$ is the pullback of the Milnor fibre $f_0^{-1}(s)$ by $\pi$ we apply to $W_s$ the procedure of the proof of  Lemma~\ref{cotanormalcrossingsglobal} and the following easy facts:
\begin{itemize}
\item The piece of Milnor fibre contained at the balls neighbouring singular points of the total transform is a union of cylinders.
\item The coverings associated to the part of Milnor fibre contained at the tubular neighbourhoods of $\dot{E}_i$ and $\dot{Z}_i$
are unramified.
\item Each set $\dot{Z}_i$ is a puntured disc.
\end{itemize}

\textsc{General case.} We reduce the proof to the previous particular case by a deformation argument. Recall that
$\tau_v$ denotes the translation in the direction of a vector $v$.
Let $v_t$ be a holomorphic family of vectors in $\CC^2$ with $v_0=O$ and such that for $t$ small enough $V(f_0\comp\tau_{v_t})$ does not
meet the origin. It is easy to check that the $2$-parameter family $F_{t,s}:=f_s\comp\tau_{v_t}$ has the following
properties:
\begin{enumerate}[(i)]
\item The set of parameters $\Delta$ such that $V(F_{t,s})$ meets the origin is a proper closed analytic subset in the parameter
space.
\item There exist positive $\eta<<\delta$ such that for any   $s$ with $0<|s|\leq \delta$ and any $t$ satisfying $0\leq |t|<\eta$ the normalization of $V(F_{t,s})\cap B_\epsilon$ is diffeomorphic to the normalization of $V(F_{0,s})=V(f_s)$.
\end{enumerate}

Choose a parametrized curve in the parameter space of the family of the form $(t(s),s)$ with $t(0)=0$ and such that for $s\neq 0$ small enough $t(s)$ is non-zero and avoids $\Delta$. Then, the normalization of $V(F_{t(s),s})$ is diffeomorphic to the normalization of $V(f_s)$ for any $s$.
Applying the particular case to the family $V(F_{t(s),s})$ we prove 
the proposition for the general case.
\end{proof}

\subsection{General global case}\label{generalglobal}
For any component $E_i$ we consider the set of irreducible components of the germ of $E_i$ at each point of $Sing((Y_0)^{red}$.
 We denote these germs by   $\{(\Gamma_k,p_k)\}_{k=1}^d$. 
We denote by $\mu_{E_i}$ the sum of Milnor numbers of these local branches, by $\nu_{E_i}$ the number of branches and define 
$$\eta_{E_i}:=\sum_{k=1}^d\sum_{l\neq k}I_{p_k}(\Gamma_k,\Gamma_l).$$
We also define the analogous numbers $\mu_{Z_i}$, $\nu_{Z_i}$ and $\eta_{Z_i}$ for any divisor $Z_i$.

For any $i$ we denote by $\dot{E}_i$ (respectively $\dot{Z}_i$) the set $E_i\setminus Sing((Y_0)^{red})$ 
(respectively $Z_i\setminus Sing((Y_0)^{red})$).

\begin{prop}
\label{cotageneralglobal}
For non-zero and small enough $s$ we have
\begin{equation}
\label{precota}
\chi(\calU_s)\leq \sum_{i=0}^mc_i(\chi(\dot{Z}_i)+\theta_{Z_i})+
\sum_{i=0}^ra_i(\chi(\dot{E}_i)+\theta_{E_i})+\sum_{p\in Y_s\cap Y_0}I_p(Y_s,(Y_0)^{red}),
\end{equation}
where $\theta(Z_i)$ and $\theta(E_i)$ are defined by
$$\theta_{Z_i}:=\nu_{Z_i}-\mu_{Z_i}-\eta_{E_i}-Z_i\centerdot((Y_0)^{red}-Z_i),$$
$$\theta_{E_i}:=\nu_{E_i}-\mu_{E_i}-\eta_{Z_i}-E_i\centerdot((Y_0)^{red}-E_i).$$
\end{prop}
\begin{proof}
The proof follows the scheme of the proof of Lemma~\ref{cotanormalcrossingsglobal}. We consider small Milnor balls
around the singular points of $(Y_0)^{red}$ and small tubular neighbourhoods around the connected components of
the complement of these balls in $(Y_0)^{red}$. We split $\calU_s$ into pieces, each being the part that maps into one of
the neighbourhoods just defined. We bound the Euler characteristic of the parts corresponding to tubular neighbourhoods
using Hurwitz formula as in the proof of Lemma~\ref{cotanormalcrossingsglobal}. We bound the Euler characteristic of 
the pieces corresponding to the Milnor balls using Proposition~\ref{cotagenerallocal}. Summing up the contributions and rearranging terms we get the desired expression.
\end{proof}

\section{Bijectivity of the Nash map for normal surface singularities}

\begin{theo}
\label{normal}
Nash mapping is bijective for any normal surface singularity defined over an algebraically closed field of characteristic equal to $0$.
\end{theo}
\begin{proof}
The argument in paragraph~\ref{lefschetz} shows that it is enough to deal with the complex case.

Let $(X,O)$ be a complex normal surface singularity. If Nash mapping is not bijective then, by Theorem~\ref{nashwedges} 
there exists a wedge $\alpha$ realizing an adjacency from
a component $E_j$ of the exceptional divisor of the minimal resolution to a different component $E_0$. We take a good representative $\alpha|_{\calU}$ and define the divisors $Y_0$ and $Y_s$ as in paragraph~\ref{familiadivisores}. As explained there, since $\calU_s$ is a disc, the lifting 
$$\tilde{\alpha}_s:\calU_s\to \tilde{X}$$ 
is the normalization of $Y_s$. We will use estimates of Section~\ref{sec:eulerchar} to get a contradiction with the fact that the Euler characteristic of $\calU_s$ is 1. In this way we show the non-existence of $\alpha$ and, by 
Theorem~\ref{nashwedges}, that the Nash mapping is bijective for normal surface singularities.

\subsection{}\label{eulercharmod}
We are going to improve slightly the estimate for $\chi(\calU_s)$ given in Proposition~\ref{cotageneralglobal}.
Remember that since $\alpha$ realizes an adjacency from $E_j$ to $E_0$ with $j\neq 0$, the divisor $E_0$ appears in $Y_0$ (see Corollary \ref{a_0}). Besides, we have a single $Z_i$, that is $Z_0$ in (\ref{decomposition}), which has the topology of a disc and intersects 
transversely $E_0$ at a smooth point of $E$. Moreover the divisor $Y_0$ is reduced at the generic point
of $Z_0$. We split $\calU_s$ in two pieces as follows. Let $\tilde{X}_1$ be a small compact tubular neighbourhood 
around the disc $Z_0$ in $\tilde{X}$. Define $\tilde{X}_2$ as the closure of the complement of $\tilde{X}_1$ in 
$\tilde{X}$. For $s$ non-zero and small enough the divisor $Y_s$ meets transversely the boundaries of $\tilde{X}_1$ and
$\tilde{X}_2$. For $i=1,2$ define $\calU_s^i$ as the normalization of $Y_s\cap\tilde{X}_i$. Since the intersection
$\calU_s^1\cap \calU_s^2$ is a union of circles we have that
\begin{equation}
\label{suma_cotas}\chi(\calU_s)=\chi(\calU_s^1)+\chi(\calU_s^2).
\end{equation}
Let us give an improved bound for the Euler characteristic of $\calU_s^1$ using the methods of paragraph ~\ref{normalcrossingslocal}. 
Let $B$ be a Milnor ball for $Y_0$ around the point $p=E_0\cap Z_0$. We may choose local coordinates $(x,y)$ around 
$p$ so that we have $E_0=V(y)$ and $Z_0=V(x)$. 
Let $g_s$ be the family of functions defining the divisor $Y_s$ locally around $p$. We have, up to a unit, the equality  
$g_0=xy^{a_0}$. The Euler characteristic of $\calU_s^1$ is bounded by the number of topological discs in the normalization of $V(g_t)\cap B$. In principle the number of circles in $\partial B\cap Y_t$ is at most $a_0+1$. There certainly appears one
circle $K_s$ which is a small deformation of $V(x)\cap\partial B$. By the connectivity of $\calU_s$, the boundary of the 
connected component of $\calU_s^1$ containing $K_s$ can not consist only of $K_s$. This implies that the maximal number of discs that can appear in $\calU_s^1$ is $a_0-1$ and hence
\begin{equation}
\label{cota1}
\chi(\calU_s^1)\leq a_0-1.
\end{equation}

The Euler characteristic of $\calU_s^2$ is bounded using Proposition~\ref{cotageneralglobal}. Notice the following identities,
$$\nu_{E_0\cap\tilde{X}_2}=\nu_{E_0}-1,$$
$$\mu_{E_0\cap \tilde{X}_2}=\mu_{E_0},$$
$$(E_0\cap\tilde{X}_2)\centerdot ((Y_0)^{red}\cap\tilde{X}_2-E_0\cap\tilde{X}_2)=E_0\centerdot ((Y_0)^{red}-E_0)-1,$$
which imply that $$\theta_{E_0\cap\tilde{X}_2}=\theta_{E_0}.$$
Then, by~(\ref{suma_cotas}) we obtain
\begin{equation}
\label{cotaA}\chi(\calU_s)\leq a_0-1+
\sum_{i=0}^ra_i(\chi(\dot{E}_i)+\theta_{E_i})+\sum_{p\in Y_s\cap Y_0\cap\tilde{X}_2}I_p(Y_s,(Y_0)^{red}).\end{equation}
Note that the last term is the total number of \emph{returns}. Defining $\delta_{a_j}=1$ if $a_j\neq 0$ and  $\delta_{a_j}=0$ if $a_j= 0$ we have the obvious bound
\begin{equation}
\label{ineq_returns}\sum_{p\in Y_s\cap Y_0\cap\tilde{X}_2}I_p(Y_s,(Y_0)^{red})\leq \sum_{j=0}^r\delta_{a_j}b_j.\end{equation}
If we denote by $k_{i,j}$ the intersection product $E_i\centerdot E_j$, by Equation~(\ref{eq:sistemalineal}) we have that  $$\sum_{j=0}^rb_j=\sum_{j=0}^r\sum_{i=0}^r\delta_{a_j}a_ik_{i,j}+1.$$ Regrouping and coming back to (\ref{ineq_returns}) we get the following 
\begin{equation}
\label{sust_ineq_returns}
\sum_{p\in Y_s\cap Y_0\cap\tilde{X}_2}I_p(Y_s,(Y_0)^{red})\leq\sum_{i=0}^ra_i(\sum_{j=0}^r\delta_{a_j}k_{i,j})+1.
\end{equation}

Now, in one hand, denoting by $g_i$ the genus of the normalization of $E_i$, we have
\begin{equation}
\label{sust_chi}
\chi(\dot{E}_i)=2-2g_i-\nu_{E_i}. 
\end{equation}
On the other hand we have that 

\begin{eqnarray}E_0\centerdot ((Y_0)^{red}-E_0)=\sum_{j\neq 0}\delta_{a_j}k_{0,j}+E_0\centerdot Z_0=\sum_{j\neq 0}\delta_{a_j}k_{0,j}+1,\nonumber\\
E_i\centerdot ((Y_0)^{red}-E_0)=\sum_{j\neq i}\delta_{a_j}k_{i,j}\ \ \ \mathrm{for\  any }\ 1\leq i\leq r.\nonumber\end{eqnarray}
and hence 
\begin{eqnarray}\label{sust_theta_0}\theta_{E_0}=\nu_{E_0}-\mu_{E_0}-\eta_{E_0}-\sum_{j\neq 0}\delta_{a_j} k_{0,j}-1,\\
\label{sust_theta_i}\theta_{E_i}=\nu_{E_i}-\mu_{E_i}-\eta_{E_i}-\sum_{j\neq i}\delta_{a_j}k_{i,j}\ \ \ \mathrm{for\  any }\ 1\leq i\leq r.
\end{eqnarray} 


Performing substitutions (\ref{sust_chi})-(\ref{sust_theta_i}) in (\ref{cotaA}) and using (\ref{sust_ineq_returns}), we get to the following:
\begin{equation}
\label{cotafinal}\chi(\calU_s)\leq \sum_{i=0}^ra_i(2-2g_i-\mu_{E_i}-\eta_{E_i}+k_{i,i}).\end{equation}

By negative definiteness, for any $0\leq i\leq r$, the self-intersection $k_{i,i}$ is a negative integer.
Observe that, since $\pi:\tilde{X}\to X$ is the minimal resolution, for any $0\leq i\leq r$, if $k_{i,i}$ is equal to
$-1$, then either the divisor $E_i$ is singular or it has positive genus 
(otherwise it is a smooth rational divisor with self-intersection equal to
$-1$ and the resolution is non-minimal). If the divisor $E_0$ has an irreducible singularity then $\mu_{E_i}$
is at least $2$. If the divisor $E_i$ has a singular point with several irreducible branches then 
$\eta_{E_i}$ is at least $2$. Therefore we have 
$$a_i(2-2g_i-\mu_i-\eta_i+k_{i,i})\leq 0$$
for any $i$ (note that $a_i\geq 0$). Hence we get that $\chi(\calU_s)\leq 0$. This is a contradiction because we know that $\calU_s$ is a disc.
\end{proof}

\section{The non-normal case}
Consider a Hironaka resolution of singularities of an algebraic variety (a resolution which is an isomorphism outside the singular locus). With any divisorial component $C$ of the exceptional locus we associate the set $N_C$ consisting of arcs in the variety
centred at the singular locus, not contained in it and whose lifting to the resolution is centred at $C$. It is an irreducible Zariski closed subset 
of the space of arcs in the variety which are not contained in the singular set. Its scheme theoretical generic point $\gamma_C$ is a formal arc in the variety
defined over the residue field of $N_C$.

Given any field $K$, a {\em formal} $K$-{\em wedge} in $X$ is a morphism
$$\alpha:Spec(K[[t,s]])\to X.$$
It can be seen as a family of arcs parametrized by $s$. Its associated special arc is the restriction of $\alpha$ to the $\{s=0\}$ locus of
$Spec(K[[t,s]])$.

We will use the following theorem: 
\begin{theo}[Theorem~5.1~\cite{Re}]\label{theo:re} An essential component $C$ of the exceptional divisor of a resolution of singularities of a variety is in the image of the Nash map
if and only if any $K$-wedge in the variety whose special arc equals the generic point $\gamma_C$ lifts to the resolution.
\end{theo}

Let $X_1$ be any reduced algebraic surface defined over a field of characteristic equal to $0$. Let
$$n:X_2\to X_1$$
be the normalization and 
$$\pi:X_3\to X_2$$ 
be the minimal resolution of the singularities of $X_2$. 

Let $\cup_{i=1}^r E_i$ be a decomposition into irreducible components of the exceptional
divisor of $\pi$. By the minimality of the resolution all these components are essential. 

Let $n^{-1}(Sing(X_1))=\cup_{i=1}^s A_i$ be a decomposition into irreducible components of 
the preimage of the singular set of $X_1$ by the normalization. Denote by $B_i$ the strict transform of $A_i$ by $\pi$. 
The decomposition into irreducible components of the exceptional divisor of the resolution $n\comp\pi$ is given by 
$$(\cup_{i=1}^s B_i)\bigcup(\cup_{i=1}^r E_i).$$
All these components are essential.

Any $K$-wedge whose special arc equals the generic point either of $N_{E_i}$ or of $N_{B_i}$ is a dominant morphism 
$$\alpha_1:Spec(K[[t,s]])\to X_1.$$
Since $Spec(K[[t,s]])$ is normal, by the universal property of the normalization, it admits a lifting
$$\alpha_2:Spec(K[[t,s]])\to X_2.$$

Assume that the special arc of $\alpha_1$ is the generic point of a certain $N_{B_j}$ for a certain $j$. Let $\eta_j$ be the generic point of $A_j$ in the scheme $X_2$. Then the image by $\alpha_2$ of the closed point of $Spec(K[[t,s]])$ is $\eta_j$. Since $\pi$ is an isomorphism in a neighbourhood of the preimage $\pi^{-1}(\eta_j)$, the wedge $\alpha_2$ lifts 
to $X_3$.

Assume that the special arc of $\alpha_1$ is the generic point of $N_{E_i}$ for a certain $i$. Then $\alpha_2$ is a wedge in $X_2$ whose special arc is the generic point
of $N_{E_i}$. Since $E_i$ is an essential component of the exceptional divisor of $\pi$, and $X_2$ is normal, case for which we have just proved that Nash mapping
is bijective, we have that $E_i$ is at the image of the Nash map for $X_2$. This implies by Theorem~\ref{theo:re} that $\alpha_2$ admits a lifting to $X_3$.

We have proved that any $K$-wedge in $X_1$ whose special arc equals the generic point of an essential component of $n\comp\pi$ lifts to $X_3$. This shows by Theorem~\ref{theo:re}
the bijectivity of the Nash map for $X_1$.

\end{document}